\documentclass[MSNbibl,number,citesort,dvips]{arxbj}
\usepackage{upgreek}
\usepackage{graphicx}

\aid{0}
\volume{20}
\issue{2}
\pubyear{2014}
\firstpage{604}
\lastpage{622}
\doi{10.3150/12-BEJ499} 

\makeatletter
\newcommand{\rright}{\right}
\newcommand{\lleft}{\left}
\newtheorem{them}{Theorem}[section]
\newtheorem{cor}[them]{Corollary}
\newcommand{\tr}{\operatorname{tr}}
\makeatother

\begin{document}
\begin{frontmatter}

\title{Information bounds for Gaussian copulas}
\runtitle{Information bounds for Gaussian copulas}

\begin{aug}
\author[1]{\fnms{Peter D.} \snm{Hoff}\corref{}\thanksref{1}\ead[label=e1]{hoff@stat.washington.edu}},
\author[2]{\fnms{Xiaoyue} \snm{Niu}\thanksref{2}}
\and
\author[1]{\fnms{Jon A.} \snm{Wellner}\thanksref{1}}
\runauthor{P.D. Hoff, X. Niu and J.A. Wellner} 
\address[1]{Departments of Statistics and Biostatistics, University of Washington, Seattle,
WA 98195-4322, USA}
\address[2]{Department of Statistics, Penn State University, University Park, PA 16802, USA}
\end{aug}

\received{\smonth{10} \syear{2011}}
\revised{\smonth{8} \syear{2012}}

%
\begin{abstract}
Often of primary interest in the analysis of multivariate data
are the copula parameters describing the dependence among the variables,
rather than the univariate marginal distributions.
Since the ranks of a multivariate dataset are invariant to
changes in the univariate marginal distributions,
rank-based estimators
are natural candidates
for semiparametric copula estimation.
Asymptotic information bounds for such estimators
can be obtained from an asymptotic analysis of the rank likelihood,
that is, the probability of the multivariate ranks.
In this article,
we obtain limiting normal distributions of
the rank likelihood for Gaussian copula models.
Our results cover models with structured correlation matrices,
such as exchangeable or circular correlation models, as
well as unstructured correlation matrices.
For all Gaussian copula models, the limiting distribution of
the rank likelihood ratio
is shown to be equal to that of a parametric likelihood ratio
for an appropriately chosen multivariate normal model.
This implies that the semiparametric information bounds
for rank-based estimators are the same as the information bounds
for estimators based on the full data, and that the multivariate normal
distributions are least favorable.
\end{abstract}

%

\end{frontmatter}

\section{Rank likelihood for copula models}

Recall that a copula is a multivariate CDF having uniform univariate
marginal distributions.
For any multivariate CDF
$F(y_1,\ldots, y_p)$ with absolutely continuous margins $F_1,\ldots, F_p$,
the corresponding copula
$C(u_1,\ldots, u_p)$ is given by
\[
C(u_1,\ldots, u_p) = F \bigl( F_1^{-1}(u_1),
\ldots, F_p^{-1}(u_p) \bigr).
\]
Sklar's theorem \cite{sklar_1959} shows that
$C$ is the unique copula for which
$F(y_1,\ldots, y_p) = C( F_1(y_1),\ldots,\allowbreak F_p(y_p))$.

In this article, we consider models consisting of
multivariate
probability distributions for which the copula is parameterized
separately from
the univariate marginal distributions. Specifically,
the models we consider consist of collections of multivariate CDFs
$\{ F(\mathbf{y}| \theta,\psi) \dvtx\mathbf{y}\in\mathbb R^p,
(\theta, \psi)
\in
\Theta\times\Psi\}$ such that
$\psi$ parameterizes the univariate marginal distributions and
$\theta$ parameterizes the copula, meaning that for a
random vector $\mathbf{Y}= (Y_1,\ldots, Y_p)^T$ with CDF $F(\mathbf
{y}|\theta,\psi)$,
\begin{eqnarray*}
\Pr(Y_j \leq y_j | \theta,\psi) & =&
F_j(y_j|\psi)\qquad\forall\theta\in\Theta, j=1, \ldots,
p,
\\
\Pr \bigl( F_1(Y_1|\psi)\leq u_1, \ldots,
F_p(Y_p|\psi) \leq u_p |\theta,\psi \bigr)
&=& C(u_1,\ldots, u_p|\theta)\qquad\forall\psi\in\Psi.
\end{eqnarray*}
We refer to such a class of distributions
as a \emph{copula-parameterized model}.
For such a model,
it will be convenient to refer to the class of copulas
$\{C(\mathbf{u}|\theta)\dvtx\theta\in\Theta\}$ as the copula model,
and the class $\{ F_1(y|\psi),\ldots, F_p(y|\psi) \dvtx\psi\in
\Psi
\}$
as the marginal model.

As an example,
the copula model for the class of $p$-variate
multivariate normal distributions
is called the Gaussian copula model, and is parameterized by letting
$\Theta$ be the set of $p\times p$ correlation matrices.
The marginal model for the $p$-variate normal distributions is the set
of all $p$-tuples of univariate
normal distributions.
The copula-parameterized models we focus on in this article are
semiparametric Gaussian copula models \cite{klaassen_wellner_1997},
for which the copula model is Gaussian and the marginal model
consists of the
set of all $p$-tuples of absolutely continuous univariate CDFs.

Let $\mathbf{Y}$ be an $n\times p$ random matrix whose rows $\mathbf
{Y}_1,\ldots,
\mathbf{Y}_n$
are i.i.d. samples
from a $p$-variate population.
We define the multivariate rank function $R(\mathbf{Y})\dvtx\mathbb
R^{n\times p}
\rightarrow\mathbb R^{n\times p}$ so that $R_{i,j}$, the $(i,j)$th
element of $R(\mathbf{Y})$, is the rank
of $Y_{i,j}$ among $\{ Y_{1,j},\ldots, Y_{n,j} \}$.
Note that the ranks $R(\mathbf{Y})$ are invariant to strictly increasing
transformations of the columns of $\mathbf{Y}$, and therefore the probability
distribution of $R(\mathbf{Y})$ does not depend on the univariate
marginal distributions
of the $p$ variables. As a result,
for any copula parameterized model and data matrix $\mathbf{y} \in
\mathbb
R^{n\times p}$
with ranks $R(\mathbf{y})=\mathbf{r}$,
the likelihood $L(\theta, \psi\dvtx\mathbf{y})$ can be decomposed as
%
\begin{eqnarray}
\label{eq:rld} L(\theta, \psi\dvtx\mathbf{y}) &=& p(\mathbf{y} | \theta,\psi) =
\Pr \bigl( R(\mathbf{Y}) = \mathbf{r} | \theta, \psi \bigr) \times p(\mathbf {y} |
\theta,\psi,\mathbf{r} )
\nonumber
\\[-8pt]
\\[-8pt]
& \equiv& L( \theta\dvtx\mathbf{r} ) \times L \bigl( \theta,\psi \dvtx[
\mathbf{y} | \mathbf{r}] \bigr) ,
\nonumber
\end{eqnarray}
where $p(\mathbf{y}|\theta, \psi)$ is the joint density of $\mathbf
{Y}$ and
$p(\mathbf{y}|\theta, \psi,\mathbf{r})$ is the conditional density
of $\mathbf{Y}$
given $R(\mathbf{Y}) = \mathbf{r}$. The function $L(\theta\dvtx
\mathbf{r}) = \Pr(R(\mathbf{Y})=\mathbf{r}|\theta)$
is called
the \emph{rank likelihood function}.
In situations where $\theta$ is the parameter of interest and $\psi$ a
nuisance parameter,
inference for $\theta$ can be obtained from the
rank likelihood function without having to estimate the margins
or specify a marginal model.
A univariate rank likelihood function was proposed by
Pettitt \cite{pettitt_1982} for estimation in monotonically transformed
regression models. Asymptotic properties of the rank likelihood for
this regression model
were studied by Bickel and Ritov \cite{bickel_ritov_1997},
and a parameter estimation scheme based on Gibbs sampling was
provided in \cite{hoff_2008b}. Rank likelihood estimation of copula
parameters was studied in
\cite{hoff_2007a}, who also extended the rank likelihood to
accommodate multivariate data with mixed continuous and
discrete marginal distributions.

The rank likelihood is constructed from the marginal probability of
the ranks and can therefore be viewed as a type of marginal likelihood.
Marginal likelihood procedures are often used for estimation in the presence
of nuisance parameters (see Section~8.3 of \cite{severini_2000} for a review).
Ideally, the statistic that generates a marginal likelihood is ``partially
sufficient'' in the sense that it contains all of the information
about the parameter of interest that can be quantified without specifying
the nuisance parameter. Notions of partial sufficiency include
$G$-sufficiency \cite{barnard_1963} and $L$-sufficiency
\cite{remon_1984},
which are motivated by group invariance and profile likelihood, respectively.
Hoff \cite{hoff_2007a} showed that
the ranks $R(\mathbf{Y})$ are both a $G$- and $L$-sufficient statistic
in the
context of copula estimation.

Although
rank-based estimators of the copula parameter $\theta$
may be appealing for the reasons described above,
one may wonder to what extent they are efficient.
The decomposition given in (\ref{eq:rld})
indicates that rank-based estimates do not use any information about
$\theta$ contained in $L(\theta,\psi\dvtx[\mathbf{y}|\mathbf
{r}])$,
the conditional density of the data given the ranks.
For at least one copula model, this information is asymptotically negligible:
Klaassen and Wellner \cite{klaassen_wellner_1997} showed that for the
bivariate normal
copula model,
a rank-based estimator is semiparametrically efficient
and has asymptotic variance
equal to the Cram\'er--Rao information bound in the bivariate normal model,
that is, the bivariate normal model is the least favorable submodel.
Genest and Werker \cite{genest_werker_2002} studied the efficiency
properties of
pseudo-likelihood estimators for two-dimensional semiparametric copula
models and showed that
the pseudo-likelihood estimators (which are functions of the bivariate ranks)
are not in general semiparametrically efficient for non-Gaussian copulas.
Chen \textit{et al.} \cite{chen_fan_tsyrennikov_2006} proposed
estimators in general multivariate copula models
that achieve semiparametric
asymptotic efficiency but
are not based solely on the multivariate ranks.
It remains unclear whether estimators based solely on the
ranks can be asymptotically efficient in general semiparametric copula models.
In particular, it is not yet
known if maximum likelihood estimators based on rank likelihoods for
Gaussian semiparametric
copula models are semiparametrically efficient.

The potential efficiency loss of rank-based estimators
can be investigated via the limiting distribution
of an appropriately scaled rank likelihood ratio.
Generally speaking, the local asymptotic normality (LAN) of a
likelihood ratio plays an important role in the asymptotic analysis
of testing and estimation procedures. For semiparametric models,
the asymptotic variance of a LAN likelihood ratio can be related to
efficient tests
\cite{choi_hall_schick_1996} and
information bounds for
regular estimators \cite
{begun_hall_huang_wellner_1983,bickel_klaassen_ritov_wellner_1993}.
In particular, the variance of the limiting normal distribution
of a LAN rank likelihood ratio provides information bounds for
locally regular rank-based estimators of copula parameters.

In this article, we obtain the limiting normal distributions
of the rank likelihood ratio for Gaussian copula models with structured
and unstructured correlation matrices.
In the next section, we give sufficient
conditions under which the
rank likelihood is LAN. The basic result is
that the rank likelihood is LAN if there exists a good rank-measurable
approximation to a LAN submodel. For Gaussian copulas, the natural
candidate submodels are multivariate normal models, for which the
log
likelihood is quadratic in the observations. In Section~\ref{sec3}, we
identify sufficient conditions for a normal quadratic form
to have a good rank-measurable approximation. This result allows us
to identify multivariate normal submodels with likelihood
ratios that asymptotically approximate the rank likelihood ratio.
In Section~\ref{sec4}, we show that for any
smoothly parameterized Gaussian copula, the rank likelihood
ratio is LAN with an asymptotic variance equal to that of the
likelihood ratio for
the corresponding multivariate normal model with unequal marginal variances.
Since the parametric multivariate normal model is a submodel
of the semiparametric Gaussian copula model,
%
and in general the semiparametric information bound based on the full data
is higher than that of any parametric submodel, our results imply that
the bounds for rank-based estimators are equal to the semiparametric
bounds for estimators based on the full data,
and that the multivariate normal models are least favorable.
These bounds can be compared to the asymptotic variance of an estimator
to assess its asymptotic efficiency.
Via two examples, in Section~\ref{sec5} we show that
pseudo-likelihood estimators are
asymptotically
efficient for some but not all Gaussian copula models.
This is discussed further
in Section~\ref{sec6}.

\section{Approximating the rank likelihood ratio}

The
local log rank likelihood ratio is defined as
\[
\lambda_r(s) = \log\frac{ L(\theta+ s/\sqrt n \dvtx\mathbf{r} )
}{ L(\theta\dvtx\mathbf{r}) },
\]
where $L(\theta\dvtx\mathbf{r})$ is defined in (\ref{eq:rld}).
Studying $\lambda_r$ is difficult because $L(\theta\dvtx\mathbf{r})$
is the integral of a copula density over a complicated set defined
by multivariate order constraints. However,
in some cases it is possible
to obtain the asymptotic distribution of $\lambda_r$ by relating
it to the local log likelihood ratio $\lambda_y$ of an appropriate
parametric multivariate model, where
%
\begin{equation}
\label{eq:pllr} \lambda_y(s,t) = \log\frac{ L(\theta+ s/\sqrt n,
\psi+t/\sqrt n \dvtx\mathbf{y} ) }{
L(\theta,\psi\dvtx\mathbf{y}) }.
\end{equation}
This method
of identifying the asymptotic distribution of $\lambda_r$
is
analogous
to the approach taken by
Bickel and Ritov \cite{bickel_ritov_1997} in their investigation of
the rank
likelihood ratio for a univariate semiparametric regression model.

In this section,
we will show that if we can find a sufficiently good
rank-measurable approximation to $\lambda_y$, then
the limiting distribution of $\lambda_r$ will match that of $\lambda_y$.
Specifically, we prove the following theorem.
\begin{them}\label{thm:mainthm}
Let $\{ F(\mathbf{y} | \theta,\psi) \dvtx\theta\in\Theta, \psi
\in\Psi\}$
be an absolutely continuous
copula parameterized model where for given values of
$\theta$ and $s$ there exists values of $\psi$ and $t$ such that
under i.i.d. sampling from $F(\mathbf{y}|\theta,\psi)$,
\begin{enumerate}
\item$\lambda_y(s,t)$ is LAN, so that $\lambda_y(s,t) \stackrel
{d}{\rightarrow}
Z$, a normal random variable, and
\item there exists a rank-measurable approximation $\lambda_{\hat y}(s,t)$
such that $\lambda_y(s,t) - \lambda_{\hat y}(s,t) \stackrel
{p}{\rightarrow} 0 $.
\end{enumerate}
Then $\lambda_r(s) \stackrel{d}{\rightarrow} Z$ as $n\rightarrow
\infty$
under i.i.d. sampling from any population with copula $C(\mathbf
{u}|\theta)$
equal to that of $F(\mathbf{y}|\theta, \psi)$
and arbitrary absolutely continuous marginal distributions.
\end{them}

\begin{pf}
Let $L(\theta,\psi\dvtx\mathbf{y})$ be the (parametric) likelihood function
for a given dataset $\mathbf{y}\in\mathbb R^{n\times p}$.
The lack of dependence of the rank likelihood on the marginal distributions
leads to the following identity relating $\lambda_r(s)$ to
$\lambda_y(s,t)$:
%
\begin{eqnarray*}
\log\mathrm{E}_\theta \bigl[ \mathrm{e}^{\lambda_y(s,t)} | R(\mathbf{Y})
= \mathbf{r} \bigr] &=&\log\int_{R(\mathbf{y} )=\mathbf{r}} \frac{ p(\mathbf
{y} | \theta+ s/\sqrt{n}, \psi+t/\sqrt{n} ) } {
p(\mathbf{y} | \theta,\psi) }
\frac{ p(\mathbf{y}|\theta,\psi) } {
\Pr( R(\mathbf{Y}) = \mathbf{r} | \theta) } \, \mathrm{d}\mathbf{y}
\\
&=& \log\frac{ \Pr( R(\mathbf{Y}) = \mathbf{r} | \theta+s/\sqrt {n} ) }{ \Pr( R(\mathbf{Y}) =\mathbf{r} | \theta) } = \lambda_r(s).
\end{eqnarray*}
%


Now suppose we would like to describe the statistical properties of
$\lambda_r(s)$ when the matrix $\mathbf{r}$ is replaced by the ranks
$R(\mathbf{Y})$, where
the rows of $\mathbf{Y}$ are i.i.d. samples from a population with copula
$C(\mathbf{u}|\theta)$. Since the distribution of the ranks
of $\mathbf{Y}$ is invariant with respect to the univariate marginal
distributions,
the particular marginal model and values of $\psi$ and $t$
are immaterial and
can be chosen to facilitate analysis.
For each $\theta$ and $s$,
our strategy will be to
select $\psi$ and $t$
such that the replacement of $\mathbf{y}$ by a rank-measurable approximation
$\hat{\mathbf{y}}$ in Equation (\ref{eq:pllr})
results in an accurate rank-based approximation
$\lambda_{\hat y}(s,t)$
of $\lambda_y(s,t)$.
Because the resulting $\lambda_{\hat y}$ is rank-measurable, we can write
\begin{eqnarray*}
\lambda_r(s)&=& \log\mathrm{E}_{\theta} \bigl[
\mathrm{e}^{\lambda_y(s,t)} | R(\mathbf{Y}) \bigr]
\\
&=& \lambda_{\hat y}(s,t) + \log\mathrm{E}_{\theta} \bigl[
\mathrm{e}^{\lambda
_y(s,t) - \lambda_{\hat y}(s,t)} | R(\mathbf{Y}) \bigr].
\end{eqnarray*}
If the approximation of $\lambda_{y}(s,t)$ by $\lambda_{\hat y}(s,t)$
is sufficiently
accurate to make the remainder term,
$ \log\mathrm{E}_{\theta}[ \mathrm{e}^{\lambda_y(s,t) - \lambda
_{\hat
y}(s,t)} | R(\mathbf{Y})]$, converge
in probability to zero as $n\rightarrow\infty$, then the asymptotic
distribution of $\lambda_r(s)$
is determined by that of $\lambda_{\hat y}(s,t)$.
Note that $\lambda_r(s)$ does not depend on $t$,
which implies that the value of $t$
for which such an approximation is available will depend on $s$ and
$\theta$.


Let $\lambda_y$ be LAN and
$ \mathbf{Y}_1,\ldots, \mathbf{Y}_n \sim$ i.i.d. $F(\mathbf
{y}|\theta,\psi)$.
For given $s$ and $t$, we will show that if
$\lambda_y(s,t)- \lambda_{\hat y}(s,t) \stackrel{p}{\rightarrow} 0
$, then
$\log\mathrm{E}_{\theta}[ \mathrm{e}^{\lambda_y(s,t) - \lambda
_{\hat
y}(s,t) } | R(\mathbf{Y}) ]
\stackrel{p}{\rightarrow} 0$,
where here and in what follows, limits are as
$n\rightarrow\infty$ and
probabilities and expectations are calculated under $\theta$ and $\psi$
unless otherwise noted.
We note that this result was essentially proven
at the end of the proof of Theorem~1 of \cite{bickel_ritov_1997}
in the context of the regression transformation model, although
details were omitted. We include the proof here for completeness.

Let $U_n = \mathrm{e}^{\lambda_y}$, $V_n = \mathrm{e}^{\lambda
_{\hat y}}$
and $\mathbf{R}_n = R(\mathbf{Y}_1,\ldots, \mathbf{Y}_n)$, so that
the exponential of the remainder term can be written as
$\operatorname{E}[{\frac{U_n}{V_n} | \mathbf{R}_n }]$.
For any $M>1$, we can write
\begin{eqnarray*}
\biggl|\operatorname{E}\biggl[{\frac{U_n}{V_n} - 1 \Big| \mathbf{R}_n }
\biggr] \biggr| &\leq& \operatorname{E}\biggl[{\biggl| \frac{U_n}{V_n} - 1 \biggr| \Big|
\mathbf{R}_n }\biggr]
\\
&=& \operatorname{E}\biggl[{\biggl| \frac{U_n}{V_n} - 1 \biggr| 1_{({U_n}/{V_n}\leq M)} \Big| \mathbf
{R}_n }\biggr] + \operatorname{E}\biggl[{\biggl| \frac{U_n}{V_n} - 1 \biggr|
1_{({U_n}/{V_n}>M) } \Big| \mathbf{R}_n }\biggr]
\\
& \leq& \operatorname{E}\biggl[{\biggl| \frac{U_n}{V_n} - 1 \biggr| 1_{({U_n}/{V_n}\leq M)} \Big|
\mathbf{R}_n }\biggr] + \operatorname{E}\biggl[{\frac{U_n}{V_n}
1_{({U_n}/{V_n}> M)} \Big| \mathbf {R}_n }\biggr]
\\
&&{} + \operatorname{E}[{1_{({U_n}/{V_n}> M)} | \mathbf{R}_n }]
\\
& =& \operatorname{E}\biggl[{\biggl| \frac{U_n}{V_n} - 1 \biggr| 1_{({U_n}/{V_n}\leq
M)} \Big|
\mathbf {R}_n }\biggr] + V_n^{-1}
\operatorname{E}[{U_n 1_{({U_n}/{V_n}> M)} | \mathbf{R}_n }]
\\
&&{} + \Pr \biggl(\frac{U_n}{V_n} > M \Big|\mathbf{R}_n \biggr)
\\
& =& a_n+ b_n+c_n.
\end{eqnarray*}
We now show that each of $a_n$, $b_n$ and $c_n$ converge in
probability to zero. To do so, we make use of the following
facts:
\begin{enumerate}
\item$U_n/V_n = \mathrm{e}^{\lambda_y- \lambda_{\hat y}}
\stackrel{p}{\rightarrow} 1 $ by the continuous mapping theorem;
\item$U_n= \mathrm{e}^{\lambda_y}$ and $V_n^{-1}=\mathrm{e}^{ -
\lambda
_{\hat y}}$ are bounded in probability, as
$\lambda_{y}$ and $\lambda_{\hat y}$ converge in distribution;
\item$\{ U_n \dvtx n \in\mathbb N\}$ is uniformly integrable,
since $\log U_n = \lambda_y$ is LAN \cite{hall_loynes_1977a};
\item If $\operatorname{E}[{|X_n| }] \rightarrow0 $ and $Z_n$ is a random
sequence, then $\operatorname{E}[{X_n | Z_n}] \stackrel
{p}{\rightarrow} 0$.
\end{enumerate}

To see that $a_n\stackrel{p}{\rightarrow} 0$ and
$c_n\stackrel{p}{\rightarrow} 0$, note that
both
$| \frac{U_n}{V_n} - 1 | 1_{({U_n}/{V_n}\leq M)}$
and
$1_{({U_n}/{V_n}> M)}$ are bounded random variables
that converge in probability to zero, so their conditional expectations
given $\mathbf{R}_n$ converge in probability to zero as well.
For the sequence $b_n$,
note that $U_n$ is $\mathrm{O}_p(1)$ as it converges in distribution,
and $1_{({U_n}/{V_n}>M)}$ is $\mathrm{o}_p(1)$ as
$\frac{U_n}{V_n}\stackrel{p}{\rightarrow} 1$, so
$\tilde U_n = U_n 1_{({U_n}/{V_n}>M)}$ is $\mathrm{o}_p(1)$.
Now $0\leq\tilde U_n \leq U_n$ for each $n$, and $\{U_n \dvtx n\in
\mathbb N\}$
is uniformly integrable, so $\{ \tilde U_n\dvtx n\in\mathbb N\}$
is uniformly integrable
as well. This and $\tilde U_n \stackrel{p}{\rightarrow} 0 $ imply that
$\operatorname{E}[{|\tilde U_n| }] = \operatorname{E}[{\tilde U_n}]
\rightarrow0$, and so
$\operatorname{E}[{\tilde U_n |\mathbf{R}_n }] \stackrel
{p}{\rightarrow} 0$. %
Since $b_n = V_n^{-1} \operatorname{E}[{\tilde U_n|\mathbf{R}_n}]$, and
$V_n^{-1}$ is $\mathrm{O}_p(1)$, $b_n$ is $\mathrm{o}_p(1)$.

Recall our original identity
relating $\lambda_r(s)$ to $\lambda_{y}(s,t)$ and $\lambda_{\hat y}(s,t)$:
\[
\lambda_r(s) 
= \lambda_{\hat y}(s,t) + \log
\operatorname{E}\bigl[{\mathrm{e}^{\lambda
_y(s,t) - \lambda_{\hat y}(s,t)} | R(\mathbf{Y}) }\bigr].
\]
We have shown that if
$\lambda_y$ is LAN
and $\lambda_y(s,t)- \lambda_{\hat y}(s,t) \stackrel{p}{\rightarrow
} 0$
under i.i.d. sampling from $F(\mathbf{y}|\theta,\psi)$,
then the remainder term goes to zero, and so
$\lambda_{y}$,
$\lambda_{\hat y}$ and $\lambda_r$ all converge to the same normal
random variable.
If the data are being sampled from a population with the same copula
as $F(\mathbf{y}|\theta,\psi)$ but different margins, then there exists
a transformation of the data
such that $F(\mathbf{y}|\theta,\psi)$ is the distribution of the transformed
population, and the result follows.
\end{pf}


For a given copula model,
Theorem~\ref{thm:mainthm} essentially says that
the asymptotic distribution of the
log rank likelihood ratio will be the same as
that of the log likelihood ratio of
any multivariate model with the same copula, as long as
the latter admits an asymptotically accurate
rank-measurable approximation.
The task of identifying the limiting distribution of
$\lambda_r$ then becomes one of identifying a suitable
marginal model for which such an approximation to the log likelihood
ratio holds.
For multivariate normal models, the log likelihood ratio is quadratic in
the observations, and so the existence of a good rank measurable
approximation depends on the accuracy of rank-based approximations
to normal quadratic forms. In the next section, we identify a
class of quadratic forms that admit sufficiently accurate
rank-measurable approximations. In Section~\ref{sec4}, we relate these forms
to multivariate normal models for which the conditions
of Theorem~\ref{thm:mainthm} hold.

\section{Rank approximations to normal quadratic forms}\label{sec3}

Let $\mathbf{Y}_1,\ldots, \mathbf{Y}_n$ be i.i.d. random column
vectors from a
member of a class of mean-zero $p$-variate normal distributions
indexed by a correlation parameter $\theta\in\Theta$
and a variance parameter $\psi\in\Psi$.
As discussed further in the next section,
the local likelihood ratio $\lambda_y$
can be expressed as a quadratic function of
$\mathbf{Y}_1,\ldots, \mathbf{Y}_n$, taking the form
\[
\lambda_y(s,t) = \Biggl( \frac{1}{\sqrt{n}} \sum
_{i=1}^n \mathbf{Y}_i^T
\mathbf{A} \mathbf{Y}_i \Biggr) + c(\theta ,\psi,s,t) +
\mathrm{o}_p(1)
\]
for some matrix $\mathbf{A}$ which could be a function of $s$, $t$,
$\theta$
and $\psi$.
A natural rank-based approximation to $\lambda_y$ is
\[
\lambda_{\hat y}( s,t) = \Biggl( \frac{1}{\sqrt{n}} \sum
_{i=1}^n \hat{\mathbf{Y}}_i^T
\mathbf{A} \hat{\mathbf{Y}}_i \Biggr) + c(\theta,\psi,s,t),
\]
where $\{\hat Y_{i,j} \dvtx i \in\{1,\ldots, n\}, j\in\{1,\ldots,
p\}
\}$
are the (approximate) normal scores, defined by
$\mathbf{R} = R(\mathbf{Y})$ and
$\hat Y_{i,j} = \sqrt{ \operatorname{Var}[{Y_{i,j} | \psi}] }\times
\Phi^{-1}
(\frac
{R_{i,j}}{n+1} )$.
Whether or not $\lambda_{\hat y} - \lambda_y \rightarrow0$ therefore
depends on the convergence to zero
of the difference between the quadratic terms of $\lambda_{\hat y} $ and
$\lambda_y$.
In this section, we show that this difference converges to
zero under certain conditions on $\mathbf{A}$ and the covariance matrix
$\mathbf{C} = \operatorname{Cov}[{\mathbf{Y}_i}]$. Specifically, we
prove the
following theorem.
\begin{them}\label{thm:quadconv}
Let $\mathbf{Y}_1,\ldots, \mathbf{Y}_n \sim$ i.i.d.
$N_p(\boldsymbol{0}, \mathbf{C})$ where
$\mathbf{C}$ is a correlation matrix, and let $\hat Y_{i,j} = \Phi^{-1}(\frac
{R_{i,j}}{n+1})$, where $R_{i,j}$ is the rank of $Y_{i,j}$ among
$Y_{1,j},\ldots, Y_{n,j}$.
Let $\mathbf{A}$ be a matrix such that
the diagonal entries of $\mathbf{A} \mathbf{C} + \mathbf{A}^T
\mathbf{C}$ are zero.
Then
\[
\frac{1}{\sqrt{n}} \sum_{i=1}^n \bigl(
\hat{\mathbf{Y}}_i^T \mathbf{A} \hat{
\mathbf{Y}}_i - \mathbf {Y}_i^T \mathbf{A}
\mathbf{Y}_i \bigr) \stackrel{p} { \rightarrow} 0 \qquad\mbox{as $n
\rightarrow\infty$.}
\]
\end{them}
%

\begin{pf}
Let $S_n = \frac{1}{\sqrt{n}} \sum_{i=1}^n ( \hat{\mathbf{Y}}_i^T
\mathbf{A} \hat
{\mathbf{Y}}_i -
\mathbf{Y}_i^T \mathbf{A} \mathbf{Y}_i )$ and
let $\tilde{\mathbf{A}} = (\mathbf{A} +\mathbf{A}^T)/2$, so that
$\mathbf{y}^T\tilde{\mathbf{A}}\mathbf{y}=
\mathbf{y}^T \mathbf{A} \mathbf{y}$ for all $\mathbf{y} \in\mathbb
{R}^p$. Then
\begin{eqnarray*}
\hat{\mathbf{Y}}^T \mathbf{A} \hat{\mathbf{Y}} -
\mathbf{Y}^T \mathbf{A} \mathbf{Y} & = & \hat{\mathbf{Y}}^T
\tilde{\mathbf{A}} \hat{\mathbf{Y}} - \mathbf {Y}^T \tilde{
\mathbf{A}} \mathbf{Y}
\\
&=& ( \hat{\mathbf{Y}} - \mathbf{Y} )^T \tilde{\mathbf{A}} ( \hat {
\mathbf{Y}} - \mathbf{Y} ) + 2 ( \hat{\mathbf{Y}} - \mathbf{Y} )^T
\tilde{\mathbf{A}} \mathbf{Y},
\end{eqnarray*}
the latter equality holding since $\tilde{\mathbf{A}}$ is symmetric.
From this, we can write
$ S_n = Q_ n + 2 L_n $ where
\begin{eqnarray*}
Q_n &=& \frac{1}{\sqrt{n}}\sum( \hat{\mathbf{Y}}_i
- { \mathbf{Y}}_i)^T \tilde{\mathbf{A}} ( \hat{
\mathbf{Y}}_i - {\mathbf{Y}}_i),
\\
L_n &=& \frac{1}{\sqrt{n}}\sum(\hat{\mathbf{Y}}_i
- {\mathbf{Y}}_i)^T \tilde{\mathbf{A}} {
\mathbf{Y}}_i.
\end{eqnarray*}
%
We can write $Q_n$ as
%
\[
Q_n = \sum_{j=1}^p \tilde
a_{j,j} \Biggl( \frac{1}{\sqrt n} \sum_{i=1}^n
(\hat Y_{i,j} - Y_{i,j})^2 \Biggr) + \sum
_{j\neq k} \tilde a_{j,k} \Biggl( \frac{1}{\sqrt n}\sum
_{i=1}^n (\hat Y_{i,j}-
Y_{i,j}) (\hat Y_{i,k} - Y_{i,k} ) \Biggr).
\]
The squared terms converge in probability to zero by 
Theorem~1 of \cite{dewet_venter_1972},
and the cross term converges in probability to zero by the Cauchy--Schwarz
inequality.

We now find conditions on $\mathbf{A}$
under which $L_n\stackrel{p}{\rightarrow} 0 $. Note that
\[
(\hat{\mathbf{y}} - \mathbf{y} )^T \tilde{\mathbf{A}} \mathbf{y} =
\sum_{j=1}^p (\hat y_{j} -
y_{j} ) \tilde{\mathbf{a}}_j^T \mathbf{y},
\]
where $\tilde{\mathbf{a}}_1,\ldots, \tilde{\mathbf{a}}_p$ are the
rows of $\tilde
{\mathbf{A}}$.
This gives
\[
L_n = \sum_{j=1}^p
L_{n,j} \equiv\sum_{j=1}^p
\frac{1}{\sqrt{n}} \sum_{i=1}^n (\hat
Y_{i,j} - Y_{i,j} ) \tilde{\mathbf{a}}_j^T
\mathbf{Y}_i.
\]
Let $\mathbf{c}_j$ be the $j$th row of $\mathbf{C}$, the correlation
matrix of $\mathbf{Y}$.
We will show that $L_{n,j} \stackrel{p}{\rightarrow} 0$ if
$\tilde{\mathbf{a}}_j^T\mathbf{c}_j = 0$ using an argument based on
conditional
expectations.
Considering $L_{n,1}$ for example,
recall that $\operatorname{E}[{\mathbf{Y}|Y_1}] = \mathbf{c}_1 Y_1$
and so
\begin{eqnarray*}
\operatorname{E}[{L_{n,1} | Y_{1,1},\ldots, Y_{n,1}
}] & = & \frac{1}{\sqrt{n}} \sum(\hat Y_{i,1} - Y_{i,1}
) \operatorname {E}\bigl[{\tilde{\mathbf{a}}_1^T
\mathbf{Y}_i | Y_{i,1} }\bigr]
\\
&=& \frac{1}{\sqrt{n}} \sum(\hat Y_{i,1} - Y_{i,1} )
\tilde{\mathbf{a}}_1^T \mathbf{c}_1
Y_{i,1} =0
\end{eqnarray*}
if $\tilde{\mathbf{a}}_j^T \mathbf{c}_j=0$. The conditional
expectation of
$L_{n,1}^2$ is given by
\begin{eqnarray*}
&&\operatorname{E}\bigl[{L_{n,1}^2 | Y_{1,1},\ldots,
Y_{n,1} }\bigr]
\\
&&\quad=\frac
{1}{n} \sum
_{i=1}^n ( \hat Y_{i,1} - Y_{i,1}
)^2 \operatorname{E}\bigl[{ \bigl(\tilde{\mathbf {a}}_1^T
\mathbf{Y}_i \bigr)^2 | Y_{i,1} }\bigr]
\\
&&\qquad{}+ \frac{1}{n} \mathop{\sum\sum}_{i_1\neq i_2} (\hat
Y_{i_1,1}-Y_{i_1,1}) (\hat Y_{i_2,1}-Y_{i_2,1})
\operatorname {E}\bigl[{\tilde{\mathbf{a}}_1^T
\mathbf{Y}_{i_1} | Y_{i_1,1} }\bigr] \operatorname{E}\bigl[{
\tilde{\mathbf{a}}_1^T\mathbf{Y}_{i_2} |
Y_{i_2,1} }\bigr].
\end{eqnarray*}
The expectations in the second sum are both proportional to $\tilde
{\mathbf{a}}_1^T\mathbf{c}_1=0$, leaving
\[
\operatorname{Var}[{L_{n,1} | Y_{1,1},\ldots, Y_{n,1}
}] = \operatorname{E}\bigl[{L_{n,1}^2 | Y_{1,1},
\ldots, Y_{n,1} }\bigr] = \frac{1 }{n} \sum
_{i=1}^n (\hat Y_{i,1}-Y_{i,1})^2
\operatorname{E}\bigl[{ \bigl(\tilde{\mathbf{a}}_1^T
\mathbf{Y}_i \bigr)^2| Y_{i,1} }\bigr].
\]
The conditional expectation
$\operatorname{E}[{(\tilde{\mathbf{a}}_1^T \mathbf{Y}_i )^2|
Y_{i,1} }]$ can be obtained
by noting that if $\mathbf{Y}\sim N_p(\boldsymbol{0}, \mathbf{C})$, then
the conditional distribution of $\mathbf{Y}$ given $Y_1$ can be
expressed as
\[
\mathbf{Y} | Y_1 \stackrel{d} {=} \mathbf{c}_1
Y_1 + \mathbf{G} \boldsymbol{\varepsilon},
\]
where $\mathbf{G}\mathbf{G}^T = \mathbf{C} - \mathbf{c}_1\mathbf
{c}_1^T$ and $\boldsymbol{\varepsilon}$ is $p$-variate
standard normal. The desired second moment is then
\begin{eqnarray*}
\operatorname{E}\bigl[{ \bigl(\tilde{\mathbf{a}}_1^T
\mathbf{Y} \bigr)^2| Y_1 }\bigr] &=& \tilde{
\mathbf{a}}_1^T \operatorname{E}\bigl[{\mathbf{Y}
\mathbf{Y}^T | Y_1 }\bigr] \tilde{\mathbf{a}}_1
\\
&=& \tilde{\mathbf{a}}_1^T \operatorname{E}
\bigl[{Y_1^2 \mathbf{c}_1
\mathbf{c}_1^T + 2 Y_1 \mathbf{c}_1
\boldsymbol {\varepsilon}^T \mathbf{G}^T + \mathbf{G}
\boldsymbol{\varepsilon }\boldsymbol {\varepsilon}^T \mathbf{G}^T
| Y_1 }\bigr] \tilde{\mathbf{a}}_1
\\
&=& \bigl( Y_1^2 -1 \bigr) \bigl(\tilde{
\mathbf{a}}_1^T\mathbf {c}_1
\bigr)^2 + \tilde{ \mathbf{a}}_1^T \mathbf{C}
\tilde{\mathbf{a}_1}
\end{eqnarray*}
which is equal to $\tilde{\mathbf{a}}_1^T\mathbf{C} \tilde{\mathbf{a}}_1$
under the condition that
$\tilde{\mathbf{a}}_1^T\mathbf{c}_1 = 0$. Letting $\gamma_1 =
\tilde{\mathbf{a}}_1^T\mathbf{C}
\tilde{\mathbf{a}}_1$, the conditional variance of $L_{n,1}$ given
the observations for the first variate is then
\[
\operatorname{Var}[{L_{n,1} | Y_{1,1},\ldots, Y_{n,1}
}] = \frac{\gamma_1 }{n}\sum(\hat Y_{i,1}-Y_{i,1})^2.
\]
Applying Chebyshev's inequality gives
\begin{eqnarray*}
\Pr\bigl( |L_{n,1}|>\varepsilon| Y_{1,1},\ldots, Y_{n,1} \bigr) &
\leq& 1 \wedge\operatorname{Var}[{L_{n,1} | Y_{1,1},\ldots,
Y_{n,1} }] /\varepsilon^2
\\
&=& 1 \wedge\frac{\gamma_1}{\varepsilon^2 } \frac{ \sum(\hat Y_{i,1} -
Y_{i,1})^2 }{n}
\\
&=& 1\wedge c_n = \tilde c_n.
\end{eqnarray*}
Now $c_n\stackrel{p}{\rightarrow}0$
as a result of Theorem~1 of \cite{dewet_venter_1972}
and therefore
so does $\tilde c_n$. But as $\tilde c_n$ is bounded,
we have $\operatorname{E}[{\tilde c_n}] \rightarrow0$, giving
\begin{eqnarray*}
\Pr\bigl( |L_{n,1}|>\varepsilon\bigr) &=& \operatorname{E}\bigl[{\Pr\bigl(
|L_{n,1}|>\varepsilon| Y_{1,1},\ldots, Y_{n,1} \bigr) }\bigr]
\\
&\leq& \operatorname{E}[{\tilde c_n }] \rightarrow0,
\end{eqnarray*}
and so $L_{n,1}\stackrel{p}{\rightarrow} 0$.
The same argument can be applied to $L_{n,j}$ for each $j$,
and so $L_n = \sum_{j=1}^p L_{n,j}\rightarrow0$ as long as $\tilde
{\mathbf{a}}_j^T\mathbf{c}_j =0$
for each $j=1,\ldots, p$, or equivalently, if the diagonal elements of
$\mathbf{A}\mathbf{C} + \mathbf{A}^T \mathbf{C}$ are zero.
\end{pf}

\section{LAN for general Gaussian copulas}\label{sec4}

In this section, we use Theorems~\ref{thm:mainthm} and~\ref{thm:quadconv}
to prove that the limiting distribution of
the rank likelihood ratio $\lambda_r$
for smoothly parameterized Gaussian copula models is
same as that
of the likelihood ratio for the corresponding normal
model with unequal marginal variances.
Specifically, we prove the following theorem.

\begin{them}\label{thm:lan}
Let $\{\mathbf{C}(\boldsymbol{\theta}) \dvtx\boldsymbol{\theta
}\in\Theta\subset\mathbb R^q\}$
be a collection
of positive definite correlation matrices such that $\mathbf
{C}(\boldsymbol{\theta})$
is twice
differentiable.
If $\mathbf{Y}_1,\ldots, \mathbf{Y}_n $
are i.i.d. from a population with absolutely continuous marginal distributions
and copula $\mathbf{C}(\boldsymbol{\theta})$ for some $\boldsymbol
{\theta}\in\Theta$, then the
distribution of the
rank likelihood ratio $\lambda_r(\mathbf{s})$ converges to a
$N(-\mathbf{s}^T
I_{\boldsymbol{\theta}\boldsymbol{\theta}\cdot\boldsymbol{\psi}}
\mathbf{s} /2, \mathbf{s}^T I_{\boldsymbol{\theta}
\boldsymbol{\theta}\cdot\boldsymbol{\psi}} \mathbf{s} )$ distribution,
where $I_{\boldsymbol{\theta}\boldsymbol{\theta}\cdot\boldsymbol
{\psi}}$ is the information for
$\boldsymbol{\theta}$
in the normal model with correlation $\mathbf{C}(\boldsymbol{\theta})$
and marginal precisions $\boldsymbol{\psi}$.
\end{them}

We note that $I_{\boldsymbol{\theta}\boldsymbol{\theta}\cdot
\boldsymbol{\psi}}$ is a function
of $\boldsymbol{\theta}$ and not of $\boldsymbol{\psi}$, as will
become clear in the proof.

\begin{pf*}{Proof of Theorem~\ref{thm:lan}}
Consider the class of mean-zero multivariate normal models with
inverse-covariance matrix $\operatorname{Var}[{\mathbf{Y}|\boldsymbol
{\theta },\boldsymbol{\psi}}]^{-1} =
\mathbf{D}(\boldsymbol{\psi})^{1/2} \mathbf{B}(\boldsymbol{\theta
}) \mathbf{D}(\boldsymbol{\psi})^{1/2}$,
where $\boldsymbol{\theta}\in\mathbb R^{q}$ and
$\mathbf{D}(\boldsymbol{\psi})$ is the diagonal matrix with diagonal elements
$\boldsymbol{\psi}\in\mathbb R^p$.
The log
probability density for a member of this class
is given by
\[
l(\mathbf{y}) = \Bigl(-p\log2\uppi+ \sum\log\psi_j + \log
|\mathbf{B}| - \mathbf{y}^T \mathbf{D}(\boldsymbol{
\psi})^{1/2} \mathbf{B} \mathbf{D}(\boldsymbol{\psi})^{1/2}
\mathbf{y} \Bigr)\big/2.
\]
The log-likelihood derivatives are
\begin{eqnarray*}
\dot{l}_{\theta_k}(\mathbf{y}) &=& \bigl[ \operatorname {tr}(
\mathbf{B}_{\theta_k} \mathbf{C} ) - \mathbf{y}^T \mathbf {D}(
\boldsymbol{\psi})^{1/2}\mathbf{B}_{\theta_k }\mathbf{D}( \boldsymbol{
\psi})^{1/2}\mathbf{y} \bigr]/2,
\\
\dot{l}_{\psi_j}(\mathbf{y}) &=& \bigl[ 1- y_j
\psi_j^{1/2}\mathbf{b}_j^T
\mathbf{D}(\boldsymbol{\psi} )^{1/2}\mathbf{y} \bigr]/(2
\psi_j),
\end{eqnarray*}
and straightforward calculations show that
\begin{eqnarray*}
I_{\boldsymbol{\psi}\boldsymbol{\psi}} &=& \mathbf{D}(\boldsymbol {\psi})^{-1} (\mathbf{I} +
\mathbf{B} \circ\mathbf{C} ) \mathbf{D}( \boldsymbol{\psi})^{-1}/4,
\\
I_{\boldsymbol{\psi}\theta_k} &=& - \mathbf{D}^{-1}(\boldsymbol {\psi})
\operatorname{diag}(\mathbf{B} \mathbf{C}_{\theta_k} )/2,
\end{eqnarray*}
where ``$\circ$'' is the Hadamard product denoting element-wise multiplication.
The local log likelihood
ratio for this model can be expressed as
\[
\lambda_y(\mathbf{s}, \mathbf{t}) = \frac{1}{\sqrt n} \sum
_{i=1}^n \mathbf{s}^T
\dot{l}_{\boldsymbol{\theta}}( \mathbf{Y}_i) + \mathbf{t}^T
\dot{l}_{\boldsymbol{\psi}}(\mathbf {Y}_i) - \frac{1}{2} \lleft[
\matrix{ \mathbf{s}
\cr
\mathbf{t}} \rright]^T I \lleft[\matrix{
\mathbf{s}
\cr
\mathbf{t}} \rright] + \mathrm{o}_p(1),
\]
which, under independent
sampling from $N_p(\boldsymbol{0}, \mathbf{D}(\boldsymbol{\psi
})^{1/2} \mathbf{C}(\boldsymbol{\theta})
\mathbf{D}(\boldsymbol{\psi})^{1/2} )$, converges in distribution to a
$N( -\mathbf{u}^TI \mathbf{u}/2, \mathbf{u}^T I \mathbf{u})$ random
variable, where
$\mathbf{u}^T = (\mathbf{s}^T, \mathbf{t}^T)$ and $I$ is the
information matrix for
$(\boldsymbol{\theta}, \boldsymbol{\psi})$.

We take
our rank based approximation $\lambda_{\hat y}$ to be equal to
$\lambda_y$
absent the $\mathrm{o}_p(1)$ term and
with each $\mathbf{Y}_i$
replaced by its approximate normal scores $\hat{\mathbf{Y}}_i$.
Clearly, we have $\lambda_{\hat y}(\mathbf{s}) -\lambda_y(\mathbf
{s}) = \mathrm
{o}_p(1)$ if
\[
\frac{1}{\sqrt{n}} \sum_{i=1}^n \bigl[
\mathbf{s}^T \dot {l}_{\boldsymbol{\theta}}(\hat{\mathbf{Y}}_i)
+ \mathbf{t}^T \dot {l}_{\boldsymbol{\psi}}(\hat{\mathbf{Y}}_i)
\bigr] - \bigl[ \mathbf{s}^T \dot{l}_{\boldsymbol{\theta}}(
\mathbf{Y}_i) + \mathbf {t}^T \dot{l}_{\boldsymbol{\psi}
}(
\mathbf{Y}_i) \bigr] = \mathrm{o}_p(1).
\]
Given $\boldsymbol{\theta}$ and $\mathbf{s}$,
we now identify a value of $\mathbf{t}$ for which the above asymptotic
result holds.
Let $\mathbf{t} = \mathbf{H} \mathbf{s}$, where $\mathbf{H} \in
\mathbb R^{p\times q}$,
so that
\begin{eqnarray*}
\mathbf{s}^T \dot{l}_{\boldsymbol\theta}(\mathbf{y}) +
\mathbf{t}^T \dot{l}_{\boldsymbol\psi}(\mathbf{y}) &=&
\mathbf{s}^T \bigl[ \dot{l}_{\boldsymbol
\theta}( \mathbf{y}) +
\mathbf{H}^T \dot{l}_{\boldsymbol\psi}(\mathbf{y}) \bigr]
\\
&=& \sum_{k=1}^q s_k \bigl[
\dot{l}_{\theta_k}(\mathbf{y}) + \mathbf{h}_k^T \dot
{l}_{\boldsymbol\psi}(\mathbf{y}) \bigr],
\end{eqnarray*}
where $\{\mathbf{h}_k, k=1,\ldots, q\}$ are the columns of $\mathbf{H}$.
Now $\dot{l}_{\theta_k}(\mathbf{y})$ and $\dot{l}_{\boldsymbol\psi
}(\mathbf{y})$
are both quadratic in $\mathbf{y}$.
Evaluating at $\boldsymbol{\psi}= \boldsymbol{1}$, we have
$\dot{l}_{\theta_k}(\mathbf{y}) = [ \operatorname{tr}(\mathbf
{B}_{\theta_k} \mathbf{C} )
- \mathbf{y}^T\mathbf{B}_{\theta_k }\mathbf{y} ]/2$ and
$\dot{l}_{\psi_j}(\mathbf{y}) = [ 1- y_j \mathbf{b}_j^T \mathbf{y}
]/2$, and so
\[
\mathbf{h}^T \dot{l}_{\boldsymbol\psi} (\mathbf{y}) = \bigl[
\mathbf{h}^T\boldsymbol{1} - \mathbf{y}^T \mathbf{D}(
\mathbf{h}) \mathbf{B} \mathbf{y} \bigr] /2.
\]
%
Therefore, we can write $ \mathbf{s}^T [ \dot{l}_{\boldsymbol\theta
}(\mathbf{y}) +
\mathbf{H}^T \dot{l}_{\boldsymbol\psi}(\mathbf{y}) ] $ as
\[
\mathbf{s}^T \bigl[ \dot{l}_{\boldsymbol\theta}(\mathbf{y}) +
\mathbf{H}^T \dot{l}_{\boldsymbol\psi}(\mathbf{y}) \bigr] = \sum
_{k=1}^q s_k \bigl[
\dot{l}_{\theta_k}(\mathbf{y}) + \mathbf{h}_k^T
\dot{l}_{\boldsymbol\psi}(\mathbf{y}) \bigr]
=\Biggl( \sum_{k=1}^q s_k
\mathbf{y}^T \mathbf{A}_k \mathbf{y} \Biggr) + c(
\mathbf{s}, \mathbf{H}, \boldsymbol{\theta}),
\]
where $c(\mathbf{s}, \mathbf{H}, \boldsymbol{\theta})$ does not
depend on $\mathbf{y}$, and
$\mathbf{A}_k$ is given by
\[
\mathbf{A}_k = - \bigl[ \mathbf{B}_{\theta_k} + \mathbf{D}(
\mathbf{h}_k) \mathbf{B} \bigr]/2.
\]
%
Substituting this representation of $\mathbf{s}^T \dot
{l}_{\boldsymbol\theta
} +
\mathbf{t}^T \dot{l}_{\boldsymbol\psi} $ into $\lambda_{\hat y}$
and $\lambda_{y}$
gives
\[
\lambda_{\hat y} - \lambda_{y} = \sum
_{k=1}^q s_k \Biggl( \frac{1}{\sqrt n}
\sum_{i=1}^n \hat{\mathbf{Y}}_i
\mathbf{A}_k \hat{\mathbf{Y}}_i - \mathbf{Y}_i
\mathbf{A}_k \mathbf{Y}_i \Biggr) +
\mathrm{o}_p(1).
\]
Theorem~\ref{thm:quadconv} implies that
this difference will converge in probability to zero if
the diagonal elements of $(\mathbf{A}_k + \mathbf{A}_k^T) \mathbf
{C}$ are zero for each
$k=1,\ldots, q$.
The value of $(\mathbf{A}_k + \mathbf{A}_k^T) \mathbf{C}$ can be
calculated as
\begin{eqnarray*}
2 \bigl(\mathbf{A}_k + \mathbf{A}_k^T
\bigr) \mathbf{C} &=& -2 \times\mathbf{B}_{\theta_k} \mathbf{C} - \mathbf{D}(
\mathbf{h}_k) \mathbf{B} \mathbf{C} - \mathbf{B} \mathbf{D}(
\mathbf{h}_k)\mathbf{C}
\\
&=& 2 \times\mathbf{B} \mathbf{C}_{\theta_k} - \bigl( \mathbf {D}(
\mathbf{h}_k) + \mathbf{B} \mathbf{D}(\mathbf{h}_k)
\mathbf{C} \bigr).
\end{eqnarray*}
The vector $\operatorname{diag}( \mathbf{D}(\mathbf{h}_k) + \mathbf
{B} \mathbf{D}(\mathbf{h}_k) \mathbf{C})$ can
be written as
\[
\operatorname{diag} \bigl( \mathbf{D}(\mathbf{h}_k) + \mathbf{B}
\mathbf{D}(\mathbf{h}_k) \mathbf{C} \bigr) = \pmatrix{
h_{k1} + \mathbf{h}_k^T (\mathbf{b}_1
\circ\mathbf{c}_1 )
\cr
\vdots\vspace*{2pt}
\cr
h_{kp} +
\mathbf{h}_k^T (\mathbf{b}_p \circ
\mathbf{c}_p )} = (\mathbf{I} + \mathbf{B}\circ\mathbf{C} )
\mathbf{h}_k,
\]
and so our condition on $\mathbf{h}_k$ becomes
\begin{eqnarray*}
(\mathbf{I}+\mathbf{B}\circ\mathbf{C} )\mathbf{h}_k &=& 2 \times
\operatorname{diag}(\mathbf{B} \mathbf{C}_{\theta_k} ),
\\
\mathbf{h}_k &=& 2 (\mathbf{I}+\mathbf{B}\circ\mathbf{C}
)^{-1} \operatorname{diag}(\mathbf{B} \mathbf{C}_{\theta_k} )
\\
&=& -I_{\boldsymbol{\psi}\boldsymbol{\psi}}^{-1} I_{\theta
_k\boldsymbol{\psi}}.
\end{eqnarray*}
Therefore, setting $\mathbf{t} =
\mathbf{H} \mathbf{s} = - I_{\boldsymbol{\psi}\boldsymbol{\psi
}}^{-1} I_{\boldsymbol{\psi}\boldsymbol{\theta}
}\mathbf{s}$ yields
a quadratic form that satisfies the conditions of Theorem~\ref{thm:quadconv}.
The result then follows via Theorem~\ref{thm:mainthm}.
The value of $\mathbf{u}^T I \mathbf{u}$ that
determines the asymptotic mean and variance of
$\lambda_y(\mathbf{s})$, $\lambda_{\hat y}(\mathbf{s})$ and
$\lambda_r(\mathbf{s})$ is given by
\begin{eqnarray*}
\mathbf{u}^T I \mathbf{u} &=& \pmatrix{ \mathbf{s}
\cr
-I_{\boldsymbol{\psi}\boldsymbol{\psi}}^{-1} I_{\boldsymbol{\psi}\boldsymbol{\theta}} \mathbf{s} }^T
\pmatrix{ I_{\boldsymbol{\theta}\boldsymbol{\theta}} & I_{\boldsymbol{\theta}\boldsymbol{\psi}}
\cr
I_{\boldsymbol{\theta}\boldsymbol{\psi}} &
I_{\boldsymbol{\psi
}\boldsymbol{\psi}} } \pmatrix{ \mathbf{s}
\cr
-I_{\boldsymbol{\psi}\boldsymbol{\psi}}^{-1}
I_{\boldsymbol{\psi
}\boldsymbol{\theta}}\mathbf{s} }
\\
&=& \mathbf{s}^T I_{\boldsymbol{\theta}\boldsymbol{\theta}} \mathbf{s} - \mathbf{s}^T
I_{\boldsymbol{\theta}\boldsymbol{\psi}} I_{\boldsymbol{\psi
}\boldsymbol{\psi}}^{-1} I_{\boldsymbol{\psi}\boldsymbol{\theta}} \mathbf{s}
\\
&=& \mathbf{s}^T I_{\boldsymbol{\theta}\boldsymbol{\theta}\cdot
\boldsymbol{\psi}} \mathbf{s}.
\end{eqnarray*}
\upqed
\end{pf*}


This result shows that the least favorable submodel of a semiparametric
Gaussian copula model is the multivariate normal model with unequal
variances, and that the information bound for any regular estimator
of $\boldsymbol{\theta}$ is given by $I_{\boldsymbol{\theta
}\boldsymbol{\theta}\cdot\boldsymbol{\psi}}$.
However, for some correlation models the value of
$I_{\boldsymbol{\theta}\boldsymbol{\theta}\cdot\boldsymbol{\psi
}}$ is equal
to the corresponding information for $\boldsymbol{\theta}$ in a model with
equal marginal variances. In such cases, the least favorable submodel
simplifies to the multivariate normal model with equal marginal variances.
To identify conditions under which this result holds,
consider the log likelihood ratio for a multivariate normal
model with equal marginal variances:
\[
\lambda_y(\mathbf{s},t) = \frac{1}{\sqrt n} \sum
_{i=1}^n \bigl[\mathbf{s}^T
\dot{l}_{\boldsymbol{\theta}}(\mathbf{Y}_i) + t \dot{l}_\psi (
\mathbf{Y}_i) \bigr] - \lleft[ \matrix{ \mathbf{s}
\cr
t }
\rright]^T I \lleft[ \matrix{ \mathbf{s}
\cr
t } \rright]\big/2 +
\mathrm{o}_p(1).
\]
Under i.i.d. sampling from $N_p(\boldsymbol{0}, \mathbf
{C}(\boldsymbol{\theta})/\psi)$, $\lambda_y(\mathbf{s},t)$
converges in distribution to a
$N(- \mathbf{u}^T I \mathbf{u}/2$, $\mathbf{u}^T I \mathbf{u})$
random variable, where
$\mathbf{u}^T = (\mathbf{s}^T, t) $ and
$I$ is the information
matrix for $(\boldsymbol{\theta}, \psi)$, for which
\begin{eqnarray*}
I_{\psi\boldsymbol{\theta}} &=& \{ I_{\psi\theta_k } \} = \bigl\{ -\operatorname{tr}(
\mathbf{B} \mathbf{C}_{\theta_k})/(2\psi) \bigr\},
\\
I_{\psi\psi} &=& p/ \bigl(2\psi^2 \bigr).
\end{eqnarray*}
%
Our candidate rank-measurable approximation to $\lambda_y(\mathbf
{s},t)$ is given
by
\[
\lambda_{\hat y}(\mathbf{s},t) = \frac{1}{\sqrt n} \sum
_{i=1}^n \bigl[\mathbf{s}^T
\dot{l}_{\boldsymbol{\theta
}}(\hat{\mathbf{Y}}_i) + t \dot{l}_\psi(
\hat{\mathbf{Y}}_i) \bigr] - \lleft[ \matrix{ \mathbf{s}
\cr
t }
\rright]^T I \lleft[ \matrix{ \mathbf{s}
\cr
t } \rright]\big/2.
\]
Recall that if for our given $\mathbf{s}$ and $\boldsymbol{\theta}$
we can find a $t$
and $\psi$ such that
$\lambda_{\hat y } - \lambda_y= \mathrm{o}_p(1)$, then the conditions
of Theorem~\ref{thm:mainthm} will be met and the asymptotic
distribution of
$\lambda_r(\mathbf{s})$ will be that of $\lambda_y(\mathbf{s},t)$.
With this in mind,
let $t = \mathbf{h}^T\mathbf{s}$ for some $\mathbf{h} \in\mathbb
R^q$, and write
$\lambda_y(\mathbf{s}, \mathbf{h}^T\mathbf{s} ) \equiv\lambda_y(\mathbf{s})$.
We will find conditions on $\mathbf{C}(\boldsymbol{\theta})$ such
that there exists
an $\mathbf{h}$ for which $\lambda_{\hat y}(\mathbf{s}) - \lambda_{y}(\mathbf{s}) =
\mathrm{o}_p(1)$,
and will show that any such $\mathbf{h}$ must be equal to
$-I^{-1}_{\psi\psi}
I_{\psi\boldsymbol{\theta}}$.
With $t=\mathbf{h}^T \mathbf{s}$ and $\psi=1$, we have
\begin{eqnarray*}
s^T \dot{l}_{\boldsymbol{\theta}}(\mathbf{y}) + t\dot{l}_\psi (
\mathbf{y}) &=& \mathbf{s}^T \bigl[ \dot{l}_{\boldsymbol{\theta}}(\mathbf{y})
+ \mathbf{h} \dot{l}_\psi(\mathbf{y}) \bigr]
\\
&=& \sum_{k=1}^q s_k \bigl[
\dot{l}_{\theta_k}(\mathbf{y}) + h_k \dot{l}_\psi(
\mathbf{y}) \bigr]
\\
&=& - \sum_{k=1}^q s_k
\mathbf{y}^T( \mathbf{B}_{\theta_k} + h_k \mathbf{B})
\mathbf{y} /2 + c(\boldsymbol{\theta}, \mathbf{s}, \mathbf{h})
\\
&=& \sum_{k=1}^q s_k
\mathbf{y}^T \mathbf{A}_k \mathbf{y} + c(\boldsymbol{
\theta}, \mathbf{s}, \mathbf{h}),
\end{eqnarray*}
where $\mathbf{A}_k = -(\mathbf{B}_{\theta_k} + h_k \mathbf{B} ) /2 =
(\mathbf{B}\mathbf{C}_{\theta_k}\mathbf{B} - h_k \mathbf{B})/2$ and
$c(\boldsymbol{\theta}, \mathbf{s}, \mathbf{h})$ does not depend on
$\mathbf{y}$.
The difference between $\lambda_{\hat y}$ and $\lambda_y$ is then
\[
\lambda_{\hat y} - \lambda_y = \sum
_{k=1}^q s_k \Biggl( \frac{1}{\sqrt{n}}
\sum_{i=1}^n \hat{\mathbf{Y}}_i
\mathbf{A}_k \hat{\mathbf{Y}}_i - {
\mathbf{Y}}_i \mathbf{A}_k {\mathbf{Y}}_i
\Biggr) + \mathrm{o}_p(1).
\]
Since $\mathbf{A}_k$ is symmetric,
Theorem~\ref{thm:quadconv}
implies that this difference will converge in probability to zero if
the diagonal elements of $\mathbf{A}_k \mathbf{C}$ are zero
for each $k=1,\ldots, q$.
This condition can equivalently be written as
follows:
\begin{eqnarray*}
\boldsymbol{0} &=& \operatorname{diag}( \mathbf{A}_k \mathbf{C} )
= \operatorname{diag}( \mathbf{B}\mathbf{C}_{\theta_k} \mathbf {B}
\mathbf{C} - h_k \mathbf{B} \mathbf{C} )/2
= \operatorname{diag}( \mathbf{B}\mathbf{C}_{\theta_k} - h_k
\mathbf{I} )/2,
\\
h_{k}\boldsymbol{1} &=& \operatorname{diag}(\mathbf{B} \mathbf
{C}_{\theta_k}).
\end{eqnarray*}
The above condition can only be met if, for each $k$,
the diagonal elements of $\mathbf{B} \mathbf{C}_{\theta_k}$
all take on a common value.
If they do, then the convergence in probability of $\lambda_{\hat
y}(\mathbf{s},t) - \lambda_y(\mathbf{s},t)$ to zero can be obtained
by setting
$t=\mathbf{h}^T\mathbf{s}$, where $h_k = \operatorname{tr}(\mathbf
{B} \mathbf{C}_{\theta_k})/p $.

Setting $\psi=1$, we have
$h_k= \operatorname{tr}(\mathbf{B} \mathbf{C}_{\theta_k})/p =
-I_{\psi\psi
}^{-1}I_{\psi\theta_k}$,
and so setting $t = \mathbf{h}^T\mathbf{s} = -I_{\psi\psi}^{-1}
I_{\psi\boldsymbol{\theta} }\mathbf{s}$
results in $\lambda_y$, $\lambda_{\hat y}$ and $\lambda_r$ each
converging in distribution to a
$N( -\mathbf{s}^TI_{\boldsymbol{\theta}\boldsymbol{\theta}\cdot
\psi}\mathbf{s}/2,\break
\mathbf{s}^TI_{\boldsymbol{\theta}\boldsymbol{\theta}\cdot\psi
}\mathbf{s})$ random variable,
where
$I_{\boldsymbol{\theta}\boldsymbol{\theta}\cdot\psi} =
I_{\boldsymbol{\theta}\boldsymbol{\theta}} - I_{\boldsymbol{\theta
}\psi} I_{\boldsymbol{\theta}\psi
}^T/I_{\psi\psi}$
is the information for $\boldsymbol{\theta}$ in this parametric model.
We summarize this result in the following corollary.
\begin{cor}\label{cor:permmat}
Let $\{\mathbf{C}(\boldsymbol{\theta}) \dvtx\boldsymbol{\theta
}\in\Theta\subset\mathbb R^q\}$
be a collection
of positive definite correlation matrices such that $\mathbf
{C}(\boldsymbol{\theta})$
is twice
differentiable, and for each $k$, the diagonal entries of
$\mathbf{B} \mathbf{C}_{\theta_k}$ are equal to some common value.
If $\mathbf{Y}_1,\ldots, \mathbf{Y}_n $
are i.i.d. from a population with absolutely continuous marginal distributions
and copula $\mathbf{C}(\boldsymbol{\theta})$ for some $\boldsymbol
{\theta}\in\Theta$, then the
distribution of the
rank likelihood ratio $\lambda_r(\mathbf{s})$ converges to a
$N(-\mathbf{s}^T
I_{\boldsymbol{\theta}\boldsymbol{\theta}\cdot\psi} \mathbf{s}
/2, \mathbf{s}^T I_{\boldsymbol{\theta}\boldsymbol{\theta}\cdot
\psi} \mathbf{s} )$ distribution,
where $I_{\boldsymbol{\theta}\boldsymbol{\theta}\cdot\psi}$ is
the information for $\boldsymbol{\theta}$
in the normal model with correlation $\mathbf{C}(\boldsymbol{\theta})$
and equal marginal precisions $\psi$.
\end{cor}

\section{Asymptotic efficiency in some simple examples}\label{sec5}

Obtaining the maximum likelihood estimator
of a copula parameter $\theta$
from the rank likelihood
is problematic due to the complicated nature of the likelihood.
An easy-to-compute alternative estimator
is the
maximizer in $\theta$
of the pseudo-likelihood, which is essentially the
probability of the observed data with the
unknown marginal CDFs replaced with
empirical estimates.
Genest \textit{et al.} \cite{genest_ghoudi_rivest_1995}
studied the asymptotic properties of this
pseudo-likelihood estimator (PLE)
and obtained a formula
for its asymptotic variance.

For Gaussian copula models, we can compare this asymptotic variance
to the information bound $I_{\theta\theta\cdot\psi}^{-1}$ obtained from
Theorem~\ref{thm:lan} to evaluate the asymptotic efficiency of the PLE.
This is most easily done in the case of a one-parameter
copula model for which the conditions of
Corollary~\ref{cor:permmat} hold, as in this case the least favorable
submodel is
a simple two-parameter multivariate normal model with equal
marginal variances.
For such models, the value of $I_{\theta\theta\cdot\psi}$ can be computed
from the variance of the efficient influence function $\check l_\theta
(\mathbf{y})$:
\begin{eqnarray*}
\check l_\theta(\mathbf{y}) &= & I_{\theta\theta\cdot\psi}^{-1}
\bigl[ \dot{l}_\theta(\mathbf{y})-I_{\theta\psi} I_{\psi\psi}^{-1}
\dot{l}_{\psi
}(\mathbf{y}) \bigr]
= I_{\theta\theta}^{-1} \bigl[\dot{l}_\theta(\mathbf{y}) -
I_{\theta\psi} \tilde l_{\psi}(\mathbf{y}) \bigr],
\end{eqnarray*}
where $\tilde l_\psi(\mathbf{y})$ is the efficient influence function for
$\psi$, given by
$ \tilde l_{\psi}(\mathbf{y}) = I_{\psi\psi\cdot\theta}^{-1} [
\dot{l}_\psi
(\mathbf{y}) -
I_{\psi\theta} I_{\theta\theta}^{-1} \dot l_{\theta}(\mathbf{y})]$
(see, e.g., \cite{bickel_klaassen_ritov_wellner_1993}, Chapter~2).
This can be compared to the influence function for the PLE,
which is given by
\[
\check l_{\theta}^P(\mathbf{y}) =
I_{\theta\theta}^{-1} \Biggl( \dot{l}_\theta(\mathbf{y}) +
\sum_{j=1}^p W_j(y_j)
\Biggr),
\]
where the likelihood derivative and information matrix are based on
the multivariate normal likelihood, and $W_j(y_j)$ is defined
as
\[
W_j(y_j) = \int_{ [0,1]^p} \biggl(
\frac{\partial^2 }{\partial\theta\,\partial u_j} \log c(\mathbf {u}|\theta) \biggr) \bigl( 1 \bigl\{
\Phi(y_j) \leq u_j \bigr\} - u_j \bigr) c(
\mathbf{u}|\theta) \,\mathrm{d}\mathbf{u}.
\]
%

By inspection,
the two influence functions are equal if
$\sum_{j=1}^p W_j(y_j) = - I_{\theta\psi} \tilde l_{ \psi}(\mathbf{y})$
$\forall\mathbf{y}\in\mathbb R^p$, in which case the
PLE is asymptotically efficient.
To compute $W_j(y_j)$ for $j=1,\ldots, p$, note that
for a Gaussian copula model, we have
\begin{eqnarray*}
\frac{\partial}{\partial\theta} \log c(\mathbf{u}|\theta) &=& - \bigl[ \tr(\mathbf{B}
\mathbf{C}_\theta) + \mathbf{y}^T \mathbf {B}_{\theta}
\mathbf{y} \bigr]/2,
\\
\frac{\partial^2}{\partial\theta\, \partial u_j }\log c(\mathbf {u}|\theta) &=& -\sum
_{k=1}^p (\mathbf{B}_{\theta})_{j,k}
\frac{ \Phi^{-1}(u_k)}{ \phi( \Phi^{-1}(u_j) ) },
\end{eqnarray*}
where $\mathbf{y} = ( \Phi^{-1}(u_1),\ldots, \Phi^{-1}(u_p) ) $,
$\mathbf{C}$
is the correlation matrix under $\theta$ and $\mathbf{B} = \mathbf{C}^{-1}$.
Straightforward calculations (\cite{shorack_2000}, page 116) give
\begin{eqnarray*}
\sum_{j=1}^p W_j(y_j)
&= & \frac{1}{2} \sum_{j=1}^p \sum
_{k=1}^p (\mathbf{B}_\theta)_{j,k}
\mathbf{C}_{j,k} \bigl(Y_j^2 -1 \bigr)
\\
&=& \tr \bigl( \mathbf{B}_\theta\mathbf{C} \bigl[\mathbf {D}(\mathbf{y}
\circ\mathbf{y} ) - \mathbf{I} \bigr] \bigr)/2 = \tr \bigl( \mathbf{B}
\mathbf{C}_\theta \bigl[ \mathbf {I} -\mathbf{D}(\mathbf{y}\circ
\mathbf{y} ) \bigr] \bigr)/2,
\end{eqnarray*}
where $\mathbf{D}(\mathbf{y} \circ\mathbf{y} )$ is the diagonal
matrix with elements
$y_1^2,\ldots, y_p^2$, and the last line follows from the
fact that $\mathbf{B}_\theta\mathbf{C} = - \mathbf{B} \mathbf
{C}_\theta$.
Recall that for the models we are considering here,
the diagonal elements of $\mathbf{B}\mathbf{C}_\theta$ are assumed
to all be equal, and
so
we can write
\[
\sum_{j=1}^p W_j(y_j)
= \frac{1}{2 p } \tr(\mathbf{B} \mathbf{C}_\theta) \sum
_{j=1}^p \bigl(1- y^2_j
\bigr).
\]
On the other hand, $I_{\theta\psi} = -\tr( \mathbf{B} \mathbf
{C}_\theta)/2$, and
so our condition for asymptotic efficiency becomes
%
\begin{eqnarray}
\label{eq:aecri} - I_{\theta\psi} \tilde l_{\psi}(\mathbf{y}) &=&
\sum_{j=1}^p W_j(y_j),
\nonumber
\\
\frac{1}{2}\tr( \mathbf{B} \mathbf{C}_\theta) \tilde
l_{\psi
}(\mathbf{y}) &=& \frac{1}{2 p } \tr(\mathbf{B}
\mathbf{C}_\theta) \sum_{j=1}^p
\bigl(1-y^2_j \bigr),
\nonumber
\\
\tilde l_{\psi}(\mathbf{y}) &=& \frac{1}{p} \sum
_{j=1}^p \bigl(1-y^2_j
\bigr).
\end{eqnarray}
We emphasize that this criterion for asymptotic efficiency only
applies to one-parameter Gaussian copula models for which the
conditions of
Corollary~\ref{cor:permmat} hold.
Such models include the one-parameter exchangeable correlation model
$\{\mathbf{C}(\theta)\dvtx\theta\in(-(p-1)^{-1},1)\}$,
for which
all off-diagonal elements are equal to $\theta$,
as well as any model in which the
rows of $\mathbf{C}(\theta)$ are permutations of one another. To see this,
note that if $\mathbf{c}_i$, the $i$th row of $\mathbf{C}(\theta)$,
is a permutation of $\mathbf{c}_j$, then $\mathbf{b}_{i}$, the $i$th
row of $\mathbf{B}$,
is the same permutation of $\mathbf{b}_j$.
Therefore $\mathbf{b}_{i}^T\mathbf{c}_{\theta,i} = \mathbf
{b}_{j}^T\mathbf{c}_{\theta,j}$ for
each $i$ and $j$,
and so the conditions of Corollary~\ref{cor:permmat} are satisfied.
Subclasses of such correlation matrices include circular correlation
models, often used for seasonal data
\cite{olkin_press_1969,khattree_naik_1994}, and any model in which the rows
of $\mathbf{C}$ are permutations of circular matrices.

\subsection*{Exchangeable correlation model}

Consider the $p=4$ exchangeable correlation matrix, for which
%
\[
\mathbf{C} = (1-\theta) \mathbf{I} + \theta\boldsymbol {1}\boldsymbol{1}^T
, \qquad\mathbf{C}_\theta= \boldsymbol{1} \boldsymbol{1}^T -
\mathbf{I} , \qquad\mathbf{B} = (1- \theta)^{-1}\mathbf{I} -
\frac{\theta}{(1-\theta)(1+3\theta)} \boldsymbol{1}\boldsymbol{1}^T.
\]
This gives
\begin{eqnarray*}
I_{\theta\theta} & =& \frac{1}{2} \operatorname{tr} (\mathbf {B}
\mathbf{C}_{\theta}\mathbf{B}\mathbf{C}_{\theta
} ) = 6
\frac{1+ 3\theta^2}{(1+ 2\theta- 3 \theta^2)^2},
\\
I_{\theta\psi} & = & -\frac{1}{2 \psi} \operatorname{tr} ( \mathbf {B}
\mathbf{C}_\theta) = \frac
{6 \theta}{1 + 2\theta- 3 \theta^2},
\\
I_{\psi\psi\cdot\theta} & = & \frac{2}{1+3\theta^2}
\end{eqnarray*}
and
\begin{eqnarray*}
\dot{l}_{\theta} (\mathbf{y}) & = & \frac{6\theta}{1+2\theta-
3\theta^2} +
\frac{1}{(1+2\theta-
3\theta^2)^2} \Biggl[ \bigl(1+3\theta^2 \bigr) \sum
_{1\le i < j \le4} y_i y_j - 3 \theta(1+\theta)
\sum_{j=1}^4 y_j^2
\Biggr],
\\
\dot{l}_{\psi} (\mathbf{y}) & = & \frac{1}{2} \biggl[
\frac{4}{\psi} - \mathbf{y}^T \mathbf{B} \mathbf{y} \biggr] = 2/
\psi- \frac{(1/2+\theta) \sum_{j=1}^4 y_j^2- \theta\sum_{1\le i
< j \le4} y_i y_j }{1 + 2\theta- 3\theta^2},
\end{eqnarray*}
so that when $\psi=1$, we have
\begin{eqnarray*}
&&\dot{l}_{\psi} (\mathbf{y}) - I_{\psi\theta} I_{\theta\theta}^{-1}
\dot{l}_{\theta}(\mathbf{y})
\\
&&\quad =  2 + \frac{\theta\sum_{1\le i <
j \le4} y_i y_j - (1/2+\theta)
\sum_{j=1}^4 y_j^2}{1+2 \theta- 3\theta^2}
\\
&&\qquad{} - \frac{6 \theta}{(1 + 2\theta-3\theta^2)} \cdot\frac
{(1+2\theta
-3\theta^2)^2}{6(1+3 \theta^2)}
\\
&& \phantom{\qquad{}-{}} {} \cdot\frac{6\theta(1 + 2\theta- 3\theta^2)
+ (1+3\theta^2) \sum_{1\le i <
j \le4} y_i y_j
- 3 \theta(1+\theta) \sum_{j=1}^4 y_j^2 }{(1+2\theta- 3\theta^2)^2}
\\
&&\quad =
\frac{2}{1+3 \theta^2} - \frac{1+2\theta-3\theta^2}{2(1+3\theta^2)(1+2\theta- 3\theta^2)} \sum_{j=1}^4
y_j^2
\\
&&\quad =  \frac{2}{1+3 \theta^2} - \frac{1}{1+3\theta^2} \frac{1}{2} \sum
_{j=1}^4 y_j^2
\\
&&\quad= \frac{1}{4} \frac{2}{1+3\theta^2} \sum_{j=1}^4
\bigl(1-y_j^2 \bigr),
\end{eqnarray*}
and so finally
\begin{eqnarray*}
\tilde{l}_{\psi} (\mathbf{y}) &=& I_{\psi\psi\cdot\theta}^{-1}
\bigl[ \dot{l}_{\psi} (\mathbf{y}) - I_{\psi\theta
}
I_{\theta\theta}^{-1} \dot{l}_{\theta} \bigr]
\\
&=& \frac{1+3\theta^2}{2} \Biggl( \frac{1}{4} \frac{2}{1+3\theta^2} \sum
_{j=1}^4 \bigl(1-y_j^2
\bigr) \Biggr)
\\
& = & \frac{1}{4} \sum_{j=1}^4
\bigl(1-y_j^2 \bigr),
\end{eqnarray*}
and so our criterion (\ref{eq:aecri}) for asymptotic efficiency is met.

%
\begin{figure}

\includegraphics{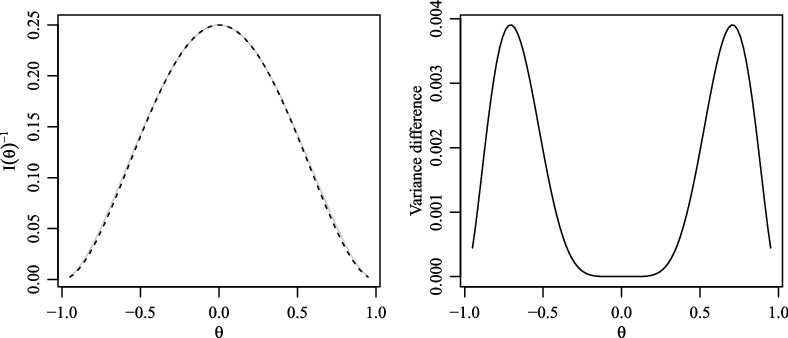}

\caption{Asymptotic variances for the circular copula model.
The left panel gives the information bound (dashed black line)
and the asymptotic variance of the PLE (gray line)
and the right panel gives the difference between these two
quantities as a function of $\theta$.}\label{fig:fig2}
\end{figure}

\subsection*{Circular correlation model}

Consider the correlation model such that
\[
\mathbf{C}(\theta) =\pmatrix{ 1 & \theta& \theta^2 & \theta
\cr
\theta& 1 & \theta& \theta^2
\cr
\theta^2 & \theta& 1 &
\theta
\cr
\theta& \theta^2 & \theta& 1 }.\vadjust{\goodbreak}
\]

For this model, we have
\[
I_{\theta\theta} = \frac{4(1+2 \theta^2)}{(1- \theta^2)^2} , \qquad I_{\theta\psi} =
\frac{4 \theta}{1 - \theta^2} , \qquad I_{\psi\psi\cdot\theta} = \frac{2}{\psi^2}
\frac{1}{1+ 2 \theta^2} .
\]
Letting $t_0 = \sum y_j^2$, $t_1=2(y_1 y_2 + y_1 y_4 + y_2y_3 + y_3y_4 )$
and $t_2=2(y_1y_3 + y_2 y_4 )$, we have
\begin{eqnarray*}
\dot{l}_\theta(\mathbf{y}) &=& \frac{4\theta}{1-\theta^2} - \frac
{ 4\theta t_0 - (1+3\theta^2)t_1 + 2\theta(1+\theta^2)t_2
}{2(1-\theta^2)^3},
\\
\dot{l}_\psi(\mathbf{y}) &=& 2/\psi- \frac{ t_0 -\theta t_1
+\theta^2 t_2 }{2(1-\theta^2)^2}.
\end{eqnarray*}
Further calculations give
\[
\tilde l_{\psi}(\mathbf{y}) = 1 - \frac{ (1-2\theta^2)t_0 + \theta^3 t_1 -\theta^2 t_2 }{4(1-\theta^2)^2 } \neq
\frac{1}{4} \sum_{j=1}^4
\bigl(1-y^2_j \bigr),
\]
and so our criterion for asymptotic efficiency is not met.
Additional calculations (available from the authors)
show that the asymptotic variance of the PLE is given by
\[
\operatorname{Var}\Biggl[{I_{\theta\theta}^{-1} \Biggl[
\dot{l}_\theta+ \sum_{j=1}^4
W_j(y_j) \Biggr] }\Biggr] = I_{\theta\theta\cdot\psi
}^{-1}
\biggl[ 1+\frac
{2\theta^6}{(1+2\theta^2)^2 } \biggr].
\]
The first panel of Figure~\ref{fig:fig2} plots the asymptotic variance
of the PLE with the
information bound, and the second panel plots their difference.
The PLE is very nearly asymptotically efficient in this example, but this
small discrepancy indicates that
the PLE is not generally asymptotically efficient for Gaussian copula
models.

\section{Discussion}\label{sec6}

In this article, we have shown that the existence of a
sufficiently accurate
rank measurable approximation to the localized log likelihood
of a
copula parameterized model
implies the local asymptotic normality of the log rank likelihood.
We have also shown that
such approximations exist for every smoothly parameterized
Gaussian copula model.
For such a copula model,
the asymptotic information bound
implied by the rank likelihood matches that
of the corresponding parametric multivariate normal submodel.
This result suggests the possibility of 
semiparametrically efficient rank-based estimators for Gaussian copula models:
Generally speaking, the information $I_r$
based on the ranks is less than or equal to the semiparametric
information $I_f$ based on the full data, as the ranks
are functions of the full data \cite{lecam_yang_1988}.
Furthermore, the semiparametric
information based on the full data is less than or equal to $I_p$, the
infimum of information
functions
over all parametric submodels, and so $I_r\leq I_f\leq I_p$ in general.
On the other hand, for Gaussian copula models
we have shown that $I_r$ is equal to
the information for a particular parametric submodel,
the corresponding multivariate normal model.
This implies that for a given Gaussian copula model, the corresponding
multivariate normal model is least favorable, that
$I_r=I_p$ and therefore $I_r=I_f=I_p$.

Based on this result, and the partial sufficiency of the multivariate
ranks in
semiparametric copula models in general,
we conjecture
that maximum likelihood estimators based on rank likelihoods are
asymptotically efficient for Gaussian copula models, and
possibly more generally whenever information bounds based on the
complete data for the
semiparametric model in question exist.
However, the rank likelihood
involves a multivariate integral over a set of
order constraints, the number of which grows with the sample size,
making it difficult to use or study.
An alternative to the rank likelihood estimator is the
pseudo-likelihood estimator
\cite{genest_ghoudi_rivest_1995},
which is a very explicit function
of the copula density, making optimization and asymptotic
analysis tractable.
For the one-parameter bivariate Gaussian
copula model, the rank-based pseudo-likelihood estimator
is asymptotically equivalent
to the normal scores correlation coefficient, which
Klaassen and Wellner \cite{klaassen_wellner_1997} showed to be
asymptotically efficient.
However,
Genest and Werker \cite{genest_werker_2002} showed with a non-Gaussian example
that the pseudo-likelihood estimator
is not generally asymptotically efficient, and in this article we have shown
that this estimator is not generally asymptotically efficient for the
restricted class of Gaussian copula models.
However, this does not rule out the possibility that
other
rank-based estimators, such as the maximizer of the rank likelihood,
are asymptotically efficient.

\section*{Acknowledgements}
Peter Hoff's research was supported in part by
NI-CHD Grant 1R01 HD067509-01A1.
Jon Wellner's research was supported in part by
 NSF Grants DMS-08-04587 and
DMS-11-04832, by NI-AID Grant 2R01 AI291968-04,
and by the Alexander von Humboldt Foundation.

%


\printhistory

\end{document}